\documentclass{amsart}
\usepackage{amsmath}
\usepackage{amscd,amssymb}
\usepackage{epsfig}

\usepackage{amsmath}
\usepackage{latexsym}
\usepackage{amssymb}
\usepackage{amsfonts}

\renewcommand{\proof}{\par\noindent{\it Proof.\ \ }}
\def\qed{\ifmmode\square\else\nolinebreak\hfill
$\Box$\fi\par\vskip12pt}

\newtheorem{thm}{Theorem}[section]

\newtheorem{lem}[thm]{Lemma}

\newtheorem{cor}[thm]{Corollary}

  \def\Rank{{\rm Rank}} \def\Res{{\rm
Res}}
\def\Ker{{\rm Ker}}

\newcommand{\Z}{\mathbb Z}
\newcommand{\C}{\mathbb C}

\newcommand{\<}{\langle}
\renewcommand{\>}{\rangle}

\begin{document}
\pagestyle{plain}

\begin{titlepage}
\begin{center}
\author{K. T. Arasu, Yu Qing Chen and Alexander Pott}
\end{center}

\address{Department of Mathematics and Statistics, Wright State
University, Dayton, Ohio 45435, USA}
\address{Institute for Algebra and Geometry,
Otto-von-Guericke-University Magdeburg, 39016 Magdeburg, Germany}
\email{}

\begin{abstract}
In this paper we prove that an abelian group contains
$(2^{2m+1}(2^{m-1}+1), 2^m(2^m+1), 2^m)$-difference sets with
$m\geqslant 3$ if and only if it contains an elementary abelian
2-group of order $2^{2m}$. Our proof shows that the method of
constructing such difference sets is essentially unique.
\end{abstract}

\title{On abelian $(2^{2m+1}(2^{m-1}+1), 2^m(2^m+1), 2^m)$-difference
sets}


\maketitle

\end{titlepage}

\section{Introduction}
Given a finite group $G$ of order $v$, a $k$-subset $D$ of $G$ is
called a $(v,k,\lambda)$-difference set if for every $g\ne 1$ in
$G$, there are exactly $\lambda$ pairs of $(d_1, d_2)\in D\times D$
such that $d_1d_2^{-1}=g$. A $(v,k,\lambda)$-difference set in a
group $G$ gives rise to a $(v,k,\lambda)$-symmetric design whose
automorphism group contains $G$ as a subgroup which acts regularly
and transitively on the design. For more details on difference sets,
we refer the reader to Beth et al. \cite{BDL}. Much of the research
on the existence or non-existence of difference sets in abelian
groups has been focused on the exponent bound of the groups.  The
best and most complete results of this type are difference sets in
abelian 2-groups \cite{Davis1, Davis2, Dillon1, Dillon2, Kraemer}.
Under the condition of character divisibility, abelian groups
containing Davis-Jedwab difference sets \cite{AC, DJ1, S} have also
been completely characterized. As for McFarland difference sets,
which are difference sets discovered by McFarland \cite{Mc} with
parameters
\begin{eqnarray*}
&&v=q^n(\frac{q^n-1}{q-1}+1),\\
&&k=q^{n-1}\frac{q^n-1}{q-1},\\
&&\lambda=q^{n-1}\frac{q^{n-1}-1}{q-1}
\end{eqnarray*}
for some positive integer $n$ and some prime power $q$, the
characterizations of abelian groups containing these difference sets
with character divisibility property are complete only for $q$ odd
\cite{MS1, MS2} or $q=4$ \cite{DJ1, MS1, MS2}.  This paper is a
contribution towards the understanding of abelian groups containing
McFarland difference sets with $q$ even and $q\geqslant 8$.
Specifically we study McFarland difference sets with $q\geqslant 8$
a power of $2$ and $n=2$. The main result of this paper is

{\thm\label{main} There exist
$(2^{2m+1}(2^{m-1}+1),2^m(2^m+1),2^m)$-difference sets in an abelian
group $G$ of order $2^{2m+1}(2^{m-1}+1)$ with $m\geqslant 3$ if and
only if $G$ contains an elementary abelian subgroup of order
$2^{2m}$.}\vspace{1mm}

The rest of the paper is organized as follows. In Section $2$, we
review the necessary background materials that are extensively used
throughout this paper. These include  characters of abelian groups,
Frattini subgroups, group algebras and building sets introduced by
Davis and Jedwab \cite{DJ1}.  The connection between building sets
and the difference sets we are interested in will also be
established. In Section $3$ we study subsets of abelian groups with
certain character sum properties. In Section $4$ we apply results
obtained in Section $3$ to study transversals of subgroups with
certain character sum properties. We introduce two types of
transversals that are important in constructing building sets for
the McFarland difference sets. Section $5$ is devoted to the proof
of Theorem~\ref{main}.

We thank Jim Davis for constructive suggestions and comments.

\section {Preliminaries}
In this section we review some of the important background materials
used in this paper including characters of abelian groups, Frattini
subgroups, group algebras and building sets.

\subsection{Characters}
A character $\chi$ of an finite abelian group $G$ is a homomorphism
from $G$ to the multiplicative group of the non-zero complex numbers
$\C^*$. The set of all characters of $G$ is a group under the
point-wise multiplication. This group is called the dual of $G$ and
is denoted by $\widehat{G}$. The identity of $\widehat{G}$, the
character $\chi$ of $G$ such that $\chi(g)=1$ for every $g\in G$, is
called the principal character of $G$ and is denoted by $\hat{1}_G$.
The group $\widehat{G}$ is isomorphic to $G$. If $H$ is a subgroup
of $G$, then restriction of a character of $G$ to $H$ is a character
of $H$ and this restriction map
$$\Res:\widehat{G}\rightarrow\widehat{H}$$ is a group
homomorphism. The kernel $\Ker(\Res)$ of the restriction map
consists of those characters $\chi$ of $G$ which are principal on
$H$, i.e. $H\leqslant\Ker(\chi)$. Hence
$$\Ker(\Res)\cong\widehat{G/H}.$$ This implies that the
restriction map $\Res$ is surjective as $G$ is finite.  The
consequence of this surjectivity is that every character $\chi$ of
$H$ can be extended to a character of $G$ and there are exactly
$G/H$ different ways of extending $\chi$ from $H$ to $G$. The
extendibility of characters of $H$ to $G$ will be frequently used in
this paper. If $G=H_1\times H_2$, then every character of $G$ is the
product of a character of $H_1$ and a character of $H_2$, i.e.
$\widehat{G}=\widehat{H_1}\times\widehat{H_2}$. To give a flavor of
how we use these facts in the subsequent sections, we prove the
following simple lemma.

{\lem\label{2elts} For any two elements $1\ne g_1\in G$ and $1\ne
g_2\in G$, there is a character $\chi\in\widehat{G}$ such that
$\chi(g_1)\ne 1$ and $\chi(g_2)\ne 1$.}\vspace{1mm}

\proof Let $K=\<g_1\>\cap\<g_2\>$. If $K\ne\{1\}$, then extension of
any non-principal characters of $K$ to $G$ gives the required
characters of $G$. If $K=\{1\}$, consider
$H=\<g_1,g_2\>\cong\<g_1\>\times\<g_2\>$. Take a non-principal
character of $\<g_1\>$ and a non-principal character of $\<g_2\>$,
their product is a character of $H$ and can be extended to $G$ to
obtain the required character.\qed

\subsection{Frattini subgroup} An element $g$ of a group $G$ is
said to be a non-generator if $G=\<X\>$ whenever $G=\<g,X\>$ for any
subset $X$ of $G$. The set of all non-generators of a group $G$ is a
subgroup. This subgroup is called Frattini subgroup and is denoted
by $\Phi(G)$. The Frattini subgroup $\Phi(G)$ of $G$ coincides with
the intersection of all proper maximum subgroups of $G$. If $G$ is a
finite $p$-group, then $\Phi(G)$ is the smallest normal subgroup $N$
of $G$ such that $G/N$ is an elementary abelian $p$-group. If $G$ is
an abelian $p$-group, then obviously $\Phi(G)$ is the homomorphic
image of the $p$-th power map of $G$. The groups $G$ we will be
investigating are $2$-groups of exponent at most $4$. The Frattini
subgroup of such a group $G$ is $\Phi(G)=\{g^2~|~g\in G\}$ and it is
elementary abelian. The characters of $G$ that are in the Frattini
subgroup of the dual of $G$ are precisely those which are principal
on the maximum elementary abelian subgroup of $G$, i.e.
$\Phi(\widehat{G})=\{\chi\in\widehat{G}~|~\chi(x)=1\mbox{ for all
}x\in G\mbox{ with }x^2=1\}$.

\subsection{Group algebras of abelian groups}

Let $G$ be an abelian group. The group algebra $\C[G]$ of the group
$G$ with complex number coefficients is the set of formal sums
$$\sum_{g\in G}a_gg,$$ where $a_g\in\C$ is the coefficient of $g$ in
the formal sum, endowed with the addition $$(\sum_{g\in
G}a_gg)+(\sum_{g\in G}b_gg)=\sum_{g\in G}(a_g+b_g)g,$$ and
multiplication $$(\sum_{g\in G}a_gg)(\sum_{g\in G}b_gg)=\sum_{g\in
G}(\sum_{h\in G}a_hb_{gh^{-1}})g.$$ It is clear that $\C[G]$ is a
$\C$-algebra. For any subset $X$ of $G$ we often identify $X$ with
the group algebra element
$$X=\sum_{g\in X}g$$ and define $X^{(t)}$ to be
$$X^{(t)}=\sum_{g\in X}g^t$$ for any integer $t\in\Z$. For
example, if $G$ is an abelian $2$-group, then $G^2=|G|G$ and
$G^{(2)}=|G/\Phi(G)|\Phi(G)$ in $\C[G]$. A $k$-element subset $D$ of
$G$ is a $(v,k,\lambda)$ difference set if and only if
$$DD^{(-1)}=(k-\lambda)+\lambda G$$ in $\C[G]$. If $f:G_1\rightarrow
G_2$ is a group homomorphism, then $f$ induces an algebra
homomorphism $f:\C[G_1]\rightarrow\C[G_2]$ by $\C$-linearly
extending $f$ from $G_1$ to $\C[G_1]$. Also if $G\cong H_1\times
H_2$ and $\chi\in\widehat{H_2}$, then $\chi$ also induces an
$\C$-algebra homomorphism $\chi:\C[G]\rightarrow\C[H_1]$. When $H_1$
is trivial, every character $\chi$ of $G$ induces a $\C$-algebra
homomorphism $\chi:\C[G]\rightarrow\C$. By the Fourier inversion
formula, two elements $X$ and $Y$ in $\C[G]$ are equal if and only
if $\chi(X)=\chi(Y)$ in $\C$ for all $\chi\in\widehat{G}$. In
particular, one has

{\lem\label{Dset} A $k$-element subset $D$ of an abelian group $G$
of order $v$ is a $(v,k,\lambda)$-difference set if and only if for
every non-principal $\chi\in\widehat{G}$,
$|\chi(D)|=\sqrt{k-\lambda}$.}\vspace{1mm}

 More generally, if
$G\cong H_1\times H_2$ and $X$ and $Y$ are in $\C[G]$, then $X=Y$ if
and only if $\chi(X)=\chi(Y)$ in $\C[H_1]$ for all
$\chi\in\widehat{H_2}$.

\subsection{Building sets}

In \cite{DJ1}, Davis and Jedwab introduced the notion of ``building
sets". This clever idea unifies the constructions of McFarland
\cite{Mc}, Spence \cite{Sp} and many others on Hadamard difference
sets, and led to the discovery of new difference sets in \cite{AC,
Chen1, Chen2, DJ1}. Given an abelian group $G$ and a subgroup $N$ of
$G$, A family of subsets $E_1$, $E_2$, $\cdots$, $E_m$ of $G$ is
called $(k,n,m)$-building sets relative to $N$ if
\begin{itemize}
\item[(i)] $|E_i|=k$, $i=1,2,\cdots, m$, \item[(ii)]
for very character $\chi$ of $G$ which is trivial on $N$ but
non-trivial on $G$, $\chi(E_i)=0$ for all $i$, and \item[(iii)]for
very character $\chi$ of $G$ which is non-trivial on $N$, there is
exactly one $i$ with $|\chi(E_i)|=n$ and $\chi(E_j)=0$ for all $j\ne
i$.
\end{itemize}
For an abelian group $G$ of order  $2^{2m+1}(2^{m-1}+1)$, the
prime $2$ is self conjugate. By Mann Test, if $D$ is a
$(2^{2m+1}(2^{m-1}+1),2^m(2^m+1),2^m)$-difference set contained in
$G$, then $2^m$ divides $\chi(D)$ for every $\chi\in\widehat{G}$,
i.e. $D$ has character divisibility property. From the results of
\cite{MS1,MS2}, The Sylow $2$-subgroup of $G$ is of exponent at
most $4$. The following lemma shows that the existence of
$(2^{2m+1}(2^{m-1}+1),2^m(2^m+1),2^m)$-difference sets in an
abelian group $G$ of order $2^{2m+1}(2^{m-1}+1)$ is equivalent to
that of $(2^{m+1}, 2^m, 2^{m-1})$-building sets in the sylow
$2$-subgroup of $G$ relative to some subgroup of order $2^m$.

{\lem\label{bset} Let $G$ be an abelian of order
$2^{2m+1}(2^{m-1}+1)$ whose Sylow $2$-subgroup has exponent no
larger than $4$. The group $G$ contains
$(2^{2m+1}(2^{m-1}+1),2^m(2^m+1),2^m)$-difference sets if and only
if its Sylow $2$-subgroup contains $(2^{m+1},2^m,2^{m-1})$-building
sets relative to some subgroup $N$ of order $2^m$.}\vspace{1mm}

\proof Let $K$ be the Sylow $2$-subgroup of $G$ and $H$ the subgroup
of $G$ of order $2^{m-1}+1$. Then $G\cong K\times H$. Suppose $D$ is
a $(2^{2m+1}(2^{m-1}+1),2^m(2^m+1),2^m)$-difference sets in $G$. The
set $D$ can be written as
$$D=E_0+h_1E_1+h_2E_2+\cdots+h_{2^{m-1}}E_{2^{m-1}},$$
where $H=\{1,h_1,h_2,\cdots,h_{2^{m-1}}\}$ and $E_0$, $E_1$,
$\cdots$, $E_{2^{m-1}}$ are subsets of $K$.  For every
$\chi\in\widehat{K}$ and every $0\leqslant i\leqslant 2^{m-1}$, by
the Fourier inversion formular, one has
\begin{eqnarray*}
\chi(E_i)&=&\frac{1}{2^{m-1}+1}\sum_{\phi\in\widehat{H}}(\chi\phi)(D)
\phi(h_i^{-1}).
\end{eqnarray*}
Since $2^m|(\chi\phi)(D)$ and $2^m$ is relatively prime to
$2^{m-1}+1$, one gets $2^m|\chi(E_i)$ for all $\chi\in\widehat{K}$
and all $0\leqslant i\leqslant 2^{m-1}$. Also when $\chi$ is
principal on $K$, the size
\begin{eqnarray*}
|E_i|&=&\frac{1}{2^{m-1}+1}\sum_{\phi\in\widehat{H}}(\hat{1}_K\phi)(D)
\phi(h_i^{-1})\\
       &\leqslant
&\frac{1}{2^{m-1}+1}\sum_{\phi\in\widehat{H}}|(\hat{1}_K\phi)(D)\phi
(h_i^{-1})|\\
       &=&\frac{1}{2^{m-1}+1}(2^m(2^m+1)+2^{m-1}2^m)\\
       &=&2^m(3-\frac{2}{2^{m-1}+1})
\end{eqnarray*}
for all $0\leqslant i\leqslant 2^{m-1}$. Since
$2^m||E_i|=\hat{1}_K(E_i)$, one has $|E_i|\leqslant 2^{m+1}$ for all
$0\leqslant i\leqslant 2^{m-1}$. The fact that $|D|=2^m(2^m+1)$
implies that $|E_i|=2^{m+1}$ for all $i$ except one, say $i=0$, with
$|E_0|=2^m$. When $\chi$ is non-principal on $K$, the absolute value
\begin{eqnarray*}
|\chi(E_i)|&=&\frac{1}{2^{m-1}+1}|\sum_{\phi\in\widehat{H}}(\chi\phi)
(D)\phi(h_i^{-1})|\\
       &\leqslant
&\frac{1}{2^{m-1}+1}\sum_{\phi\in\widehat{H}}|(\chi\phi)(D)\phi(h_i^{-
1})|\\
       &=&2^m
\end{eqnarray*}
for all $0\leqslant i\leqslant 2^{m-1}$. Since $2^m|\chi(E_i)$ and
$\chi(E_i)\in\Z[\sqrt{-1}]$ as $K$ is of exponent no larger than
$4$, one has $|\chi(E_i)|=2^m$ or $0$ for all non-principal
$\chi\in\widehat{K}$ and all  $0\leqslant i\leqslant 2^{m-1}$. This
implies that $E_0$ is a coset of a subgroup, say $N$, of order $2^m$
in $K$. From $DD^{(-1)}=2^{2m}+2^mG$, one gets
\begin{eqnarray*}
\sum_{i=0}^{2^{m-1}}E_iE_i^{(-1)}&=&2^{2m}+2^mK
\end{eqnarray*} and
\begin{eqnarray*}
\sum_{i=0}^{2^{m-1}}|\chi(E_i)|^2=\sum_{i=0}^{2^{m-1}}\chi(E_iE_i^{(-
1)})=2^{2m}
\end{eqnarray*}
for every non-principal $\chi\in\widehat{K}$. Hence $E_1$, $E_2$,
$\cdots$, $E_{2^{m-1}}$ form $(2^{m+1},2^m,2^{m-1})$-building sets
in $K$ relative to $N$. Conversely, if $K$ contains
$(2^{m+1},2^m,2^{m-1})$-building sets $E_1$, $E_2$, $\cdots$,
$E_{2^{m-1}}$ relative to $N$ of order $2^m$, then it is easy to
check that when $E_0=N$ and
$$D=E_0+h_1E_1+h_2E_2+\cdots+h_{2^{m-1}}E_{2^{m-1}},$$ one has
$|\chi(D)|=2^m$ for all non-principal $\chi\in\widehat{G}$ and
therefore, by Lemma~\ref{Dset}, $D$ is a
$(2^{2m+1}(2^{m-1}+1),2^m(2^m+1),2^m)$-difference set in $G$.\qed

\section{Subsets in abelian $2$-groups with certain character sum
properties}

In this section we study subsets $E$ in abelian $2$-groups $G$ of
exponent no larger than $4$ such that $|\chi(E)|=|E|/2$ or $0$ for
all non-principal $\chi\in\widehat{G}$. The results in this section
will used in subsequent sections in investigating
$(2^{m+1},2^m,2^{m-1})$-building sets. We start with the following
two simple technical lemmas.  The first one shows that $E\cap K$
inherits character divisibility from $E$ while the second
demonstrates that $E$ contains a coset of a subgroup of order
$|E|/4$ if $|\chi(E)|=|E|/2$ for some character $\chi$ of order $2$
in ${\widehat G}$.

{\lem\label{cong} Let $G$ be a finite abelian group of even order
and $E$ a subset of $G$. Suppose $E$ has the property that there is
a positive integer $l$ such that for every character
$\chi\in\widehat{G}$, the character sum $\chi(E)\equiv
0\;(\!\!\!\!\mod 2l)$. Then for every index $2$ subgroup $K$ of $G$,
the subset $E\cap K$ has the property that for every
$\chi\in\widehat{K}$, the character sum $\chi(E\cap K)\equiv 0\;\;
(\!\!\!\!\mod l)$.}\vspace{1mm}

\proof Let $K$ be an index $2$ subgroup of $G$ and $x\in G\setminus
K$. The set $E$ can be decomposed into $$E=(E\cap K)+x(x^{-1}E\cap
K).$$ For every character $\chi\in\widehat{K}$, there are two
characters $\chi^{+}$ and $\chi^{-}$ in $\widehat{G}$ such that
$\chi^{+}|_K=\chi^{-}|_K=\chi$ and $\chi^{+}(x)+\chi^{-}(x)=0$.
Applying $\chi^{+}$ and $\chi^{-}$ to $E$, one gets
\begin{eqnarray*}
&\chi^{+}(E)=\chi(E\cap K)+\chi^{+}(x)\chi(x^{-1}E\cap K)\equiv
0\;\:(\!\!\!\!\!\!\mod 2l),&\\
&\chi^{-}(E)=\chi(E\cap K)+\chi^{-}(x)\chi(x^{-1}E\cap K)\equiv
0\;\:(\!\!\!\!\!\!\mod 2l).&
\end{eqnarray*}
The sum of the two equations ensures that $\chi(E\cap K)\equiv
0\;(\!\!\!\!\mod l)$.\qed

{\lem\label{E2} Let $G$ be an abelian $2$-group of exponent no
larger than $4$ and $E$ a subset of $G$ such that $|E|\equiv
0\;\;(\!\!\!\!\mod 4)$ . If the subset $E$ satisfies
\begin{itemize}
\item[(1)] for every character $\chi\in\widehat{G}$, the  value
$\chi(E)\equiv 0\;\;(\!\!\!\!\mod |E|/2)$,  and
\item[(2)] there is a character $\phi\in\widehat{G}$ of order $2$ such
that $|\phi(E)|=|E|/2$,
\end{itemize}
then $E$ contains a coset of a subgroup of order $|E|/4$ in $G$ .}
\vspace{1mm}

\proof Let $K$ be the kernel of $\phi$. Since $\phi$ is of order
$2$, the subgroup $K$ is of index $2$ in $G$. Let $x\in G\setminus
K$. Then $E=(E\cap K)+x(x^{-1}E\cap K)$ and $\phi(E)=|E\cap
K|-|x^{-1}E\cap K|=\pm |E|/2$.  Hence either $|E\cap K|=|E|/4$ or
$|x^{-1}E\cap K|=|E|/4$. By Lemma~\ref{cong},  $\chi(E\cap
K)\equiv\chi(x^{-1}E\cap K)\equiv 0\;\;(\!\!\!\!\mod |E|/4)$ for all
$\chi\in\widehat{K}$. Therefore $E$ contains a coset of a subgroup
of order $|E|/4$ in $G$.\qed

The next lemma shows that the subgroup $(\Z/2\Z)^3$ in $(\Z/4\Z)^3$
has no transversal $E$ such that $\chi(E)\equiv 0\;\;(\!\!\!\!\mod
4)$ for all characters of $(\Z/4\Z)^3$.  This is actually a crucial
fact of our proof of Theorem~\ref{main}.

{\lem\label{Z43} Let $G\cong(\Z/4\Z)^3$ and $N$ the maximum
elementary abelian subgroup of $G$. There is no transversal $E$ of
$N$ in $G$ such that $\chi(E)\equiv 0\;\;(\!\!\!\!\mod 4)$ for every
$\chi\in\widehat{G}$.}\vspace{1mm}

\proof Let $G=\<x\>\times\<y\>\times\<z\>$, where $x$, $y$ and $z$
are of order 4. Up to equivalence by isomorphisms and translations
of $G$, we may assume that
$$E=1+x+y+z+xyw_{xy}+yzw_{yz}+zxw_{zx}+xyzw_{xyz},$$
where $w_{xy}$, $w_{yz}$, $w_{zx}$ and $w_{xyz}$ are elements in
$N$. Since $\chi(E)\equiv 0\;(\!\!\!\!\mod 4)$ for every
$\chi\in\widehat{G}$, one can check that $$\prod_{g\in E}\chi(g)=1$$
for every $\chi\in\widehat{G}$. Therefore
$$\prod_{g\in E}g=w_{xy}w_{yz}w_{zx}w_{xyz}=1$$ and
$w_{xyz}=w_{xy}w_{yz}w_{zx}$. If we apply character $\chi$ such that
$\chi(x)=\pm\sqrt{-1}$, $\chi(y)=\pm 1$ and $\chi(z)=\pm 1$, we find
that $1+\chi(y)+\chi(z)+\chi(y)\chi(z)\chi(w_{yz})\equiv
0\;(\!\!\!\mod 4)$ and
$1+\chi(y)\chi(w_{xy})+\chi(z)\chi(w_{zx})+\chi(y)\chi(z)\chi(w_{xyz})
\equiv
0\;(\!\!\!\!\mod 4)$. This is equivalent to $\chi(w_{yz})=1$ because
$a_1+a_2+a_3+a_4\equiv 0\;(\!\!\!\mod 4)$ if and only if
$a_1a_2a_3a_4=1$ when $a_i=\pm 1$ for $i=1,2,3,4$.  Therefore
$$w_{yz}\in\{1,y^2,z^2,y^2z^2\}.$$ Similarly,
\begin{align*}
&w_{zx}\in\{1,x^2,z^2,x^2z^2\},\\
&w_{xy}\in\{1,x^2,y^2,x^2y^2\},\\
&w_{yz}w_{xyz}\in\{1,x^2,y^2z^2,x^2y^2z^2\},\\
&w_{zx}w_{xyz}\in\{1,y^2,z^2x^2,x^2y^2z^2\},\\
&w_{xy}w_{xyz}\in\{1,z^2,x^2y^2,x^2y^2z^2\},\\
&w_{xyz}\in\{x^2,y^2,z^2,x^2y^2z^2\}.
\end{align*}
Therefore we need $w_{xy}$, $w_{yz}$, $w_{xz}$ and $w_{xyz}$ in $N$
satisfy
\begin{align}
&w_{xy}\in\{1,x^2,y^2,x^2y^2\}\cap
w_{xyz}\{1,z^2,x^2y^2,x^2y^2z^2\},\\
&w_{yz}\in\{1,y^2,z^2,y^2z^2\}\cap
w_{xyz}\{1,x^2,y^2z^2,x^2y^2z^2\},\\
&w_{zx}\in\{1,x^2,z^2,x^2z^2\}\cap
w_{xyz}\{1,y^2,z^2x^2,x^2y^2z^2\},\\
&w_{xyz}\in\{x^2,y^2,z^2,x^2y^2z^2\}\\
&w_{xy}w_{yz}w_{xz}w_{xyz}=1.
\end{align}
 One can easily check that there is
no $w_{xyz}$ in $\{x^2,y^2,z^2,x^2y^2z^2\}$ that can produce
$w_{xy}$, $w_{yz}$ and $w_{zx}$ that satisfy all conditions
above.\qed

For a $p$-group $G$ and a normal subgroup $N$ of $G$, it is clear
that $\Phi(N)$ is contained in $N\cap\Phi(G)$ as $\Phi(N)$ is the
smallest normal subgroup $X$ of $N$ such that $N/X$ is elementary
abelian while $N/(N\cap\Phi(G))$ is a subgroup of $G/\Phi(G)$ which
is elementary abelian. An interesting consequence of
Lemma~\ref{Z43}, which we prove as the following corollary, is that
if a subgroup $N$ of an abelian $2$-group $G$ of exponent no larger
than $4$ has a transversal $E$ with the property that $\chi(E)\equiv
0\;\;(\!\!\!\!\mod |E|/2)$ for all $\chi\in\widehat{G}$, then
$\Phi(N)$ is at most index $4$ in $N\cap\Phi(G)$. that is

{\cor\label{CZ43} Let $G$ be an abelian $2$-group of exponent no
larger than $4$ and $N$ a subgroup of $G$. Let $E$ be a transversal
of $N$ in $G$. If for every character $\chi\in\widehat{G}$, the
value $\chi(E)\equiv 0\;\;(\!\!\!\!\mod |E|/2)$, then
$|(N\cap\Phi(G))/\Phi(N)|\leqslant 4$.}\vspace{1mm}

\proof If $|(N\cap\Phi(G))/\Phi(N)|\geqslant 8$, there is a subgroup
$K\cong(\Z/4\Z)^3$ of $G$ such that $N\cap K\cong(\Z/2\Z)^3$ and
$\Phi(N)\cap K=\{1\}$. Then the group $N$ has a subgroup $N'$ such
that $N'\cap K=\{1\}$ and $N=(N\cap K)N'$. In the quotient group
$G/N'$, the image $\overline{E}$ of $E$ is a transversal of the
image $\overline{N}=N/N'\cong N\cap K <\overline{K}\cong K$. The
transversal also has the property that for every character
$\chi\in\widehat{\overline{G}}$, the value $\chi(\overline{E})\equiv
0\;(\!\!\!\!\mod |\overline{E}|/2)$ since $E$ has. Now we can find a
sequence of subgroups
$\overline{K}=H_0<H_1<H_2\cdots<H_{k-1}<H_k=\overline{G}$, where $k$
is the positive integer satisfies $|\overline{G}|=2^{k+6}$, such
that $|H_i/H_{i-1}|=2$ for $i=1$, $2$, $\cdots$, $k$. Applying
Lemma~\ref{cong} successively to the subgroups $H_{k-1}$, $H_{k-2}$,
$\cdots$, $H_0=\overline{K}$, one finds that
$\overline{E}\cap\overline{K}$ has the property that
$\chi(\overline{E}\cap\overline{K})\equiv 0$ $(\!\!\!\!\mod 4)$ for
every character $\chi\in\widehat{\overline{K}}$. Also
$\overline{E}\cap\overline{K}$ is a transversal of $\overline{N}$ in
$\overline{K}$. This contradicts Lemma~\ref{Z43}.\qed

As an application of Corollary~\ref{CZ43}, we prove the following
lemma which shows that every set in $(2^{m+1},2^m,2^{m-1})$-building
sets satisfies the conditions of Lemma~\ref{E2}.

{\lem\label{Eg} Let $m\geqslant 3$ be an integer. Let $G$ be an
abelian group of order $2^{2m+1}$ and of exponent at most $4$, and
$E$ a subset of $G$ of size $2^{m+1}$. If for every non-principal
character $\chi\in\widehat{G}$, $|\chi(E)|=2^m$ or $0$, then there
is a character $\chi$ of order $2$ in $\widehat{G}$ such that
$|\chi(E)|=|E|/2$ .}\vspace{1mm}

\proof Suppose to the contrary. Then $\chi(E)=0$ for every character
$\chi$ of order $2$ in $\widehat{G}$, i.e. every
$\chi\in\widehat{G}$ with $\Phi(G)\leqslant\Ker(\chi)$. Let
$\overline{G}=G/\Phi(G)$. Then the image $\overline{E}$ of $E$ in
$\overline{G}$ must be a multiple of the whole group. Hence
$|\overline{G}|=|G/\Phi(G)|\leqslant |E|=2^{m+1}$, or equivalently,
$|\Phi(G)|\geqslant |G|/|E|=2^m$. Since $G$ is of exponent at most
$4$, $|\Phi(G)|\leqslant 2^m$. Therefore $|\Phi(G)|=2^m$ and $E$ is
a transversal of $\Phi(G)$ in $G$, which contradicts
Corollary~\ref{CZ43}.\qed

\section{Transversals in abelian $2$-groups with certain character sum
properties}

In this section we introduce two types of transversals and show that
they arise frequently in the building sets that we are interested
in. Let $G$ be an abelian $2$-group of exponent no larger than $4$
and $N$ a subgroup of $G$. Let $E$ be a transversal of $N$ in $G$.
The two types of transversals that we want to investigate are
\begin{itemize}
\item[(I)]There are two subgroups $H_1$ and $H_2$ of $G$ such that
$|H_1|=|H_2|=|E|/2$ and $E=aH_1+bH_2$ for some $a$ and $b$ in $G$;
\item[(II)]There is a subgroup $H$ of $G$ of order $|E|/8$ and an
$8$-element subset $E{'}$ of $G$ such that $E=HE{'}$.
\end{itemize}

The next two two lemmas give descriptions of the structure of type
(I) transversals. The first one is valid without any restrictions
while the second is for transversal $E$ with the property that
$|\chi(E)|=|E|/2$ or $0$.

{\lem\label{I1} If $E$ is of type $\rm{(I)}$, then $|H_1\cap
N|=|H_2\cap N|=1$, $H_1N=H_2N$ and $a^{-1}b\notin
H_1N$.}\vspace{1mm}

\proof Since $E$ is a transversal of $N$ in $G$, it is obvious that
$|H_1\cap N|=|H_2\cap N|=1$. Since $a^{-1}E$ is also a transversal
of $N$ in $G$, the intersection $NH_1\cap a^{-1}E=H_1$ is a
transversal of $N$ in $NH_1$, and therefore $a^{-1}bH_2\cap
NH_1=\emptyset$. Hence $a^{-1}b\notin NH_1$. Since $NH_1$ is an
index $2$ subgroup of $G=NE$, one finds that the subgroup $H_2
\subseteq ab^{-1}(G\setminus NH_1)=NH_1$. Since $|H_1N|=|H_2N|$, one
gets $H_1N=H_2N$.\qed

{\lem\label{I} Let $E$ be a transversal of type $\rm{(I)}$. If $E$
has the property that $|\chi(E)|=|E|/2$ or $0$ for every
non-principal character $\chi\in\widehat{G}$, then
$N\leqslant\<H_1,H_2\>$, $|H_1\cap H_2|=|E|/(2|N|)$ and
$H_1/(H_1\cap H_2)\cong H_2/(H_1\cap H_2)\cong N$.}\vspace{1mm}

\proof Let $H=\<H_1,H_2\>$, by Lemma~\ref{I1}, $H\leqslant
NH_1=NH_2$. From the absolute values of character sums of $a^{-1}E$,
we know that $G=\<a^{-1}E\>=\<a^{-1}b,H\>$. We show that $N\leqslant
H$. If there is an $x\in N$ such that $x\notin H$, then
$x=(a^{-1}b)^sh$ for some integer $1\leqslant s\leqslant 3$ and some
$h\in H$. By Lemma~\ref{I1}, $s\ne 1$ or $3$ as
$(a^{-1}b)^s=xh^{-1}\in NH=NH_1=NH_2$. This  implies that $s=2$,
$a^{-1}b$ is an element of order $4$ and $|\<a^{-1}b\>\cap H|=1$.
Hence there is a non-principal character $\chi\in\widehat{G}$ with
$\chi((a^{-1}b)^2)=-1$ and $\chi$ is principal on $H$, so that
$|\chi(E)|=\sqrt{2}|E|/2$. This contradicts the assumption that
$|\chi(E)|=|E|/2$ or $0$ for every non-principal character
$\chi\in\widehat{G}$. Therefore $N\leqslant H$ and $H=NH_1=NH_2$.
>From $|\<H_1,H_2\>|=|H_1||H_2|/|H_1\cap H_2|$, one gets $|H_1\cap
H_2|=|E|/(2|N|)$. Now consider the group $\overline{G}=G/(H_1\cap
H_2)$ and the image $\overline{E}$ of $E$ in $\overline{G}$.
Clearly, $\overline{E}=|H_1\cap
H_2|(\overline{a}\overline{H}_1+\overline{b}\overline{H}_2)$ and
$\overline{a}\overline{H}_1+\overline{b}\overline{H}_2$ is a
transversal of $\overline{N}$ in $\overline{G}$ which has the
property that for every non-principal character
$\chi\in\widehat{\overline{G}}$, the absolute value
$|\chi(\overline{a}\overline{H}_1+\overline{b}\overline{H}_2)|=
|\overline{a}\overline{H}_1+\overline{b}\overline{H}_2|/2$ or $0$,
where $\overline{H}_1=H_1/(H_1\cap H_2)$,
$\overline{H}_2=H_2/(H_1\cap H_2)$ and $\overline{N}\cong N$ is the
image of $N$ in $\overline{G}$. Since
$|\overline{H}_1\cap\overline{H}_2|=|\overline{H}_1\cap\overline{N}|
=|\overline{H}_2\cap\overline{N}|=1$ and
$\overline{H}_1\overline{H}_2=\overline{H}_1\overline{N}
=\overline{H}_2\overline{N}$, one gets
$\overline{H}_1\cong\overline{H}_2\cong\overline{N}\cong N$. \qed

The next two lemmas show that for a type (II) transversal $E$ of $N$
in $G$, there is a restriction on the size of the subgroup $N$ when
we impose the condition $\chi(E)=|E|/2$ or $0$ for every
non-principal $\chi\in\widehat G$.

{\lem\label{size} Let $G$ be a finite abelian $2$-group and $N$ a
subgroup of $G$. If a transversal $E$ of $N$ in $G$ satisfies
\begin{itemize}
\item[(1)] $|E|\leqslant 8$, and
\item[(2)] $|\chi(E)|=|E|/2$ or $0$ for every non-principal character
$\chi\in\widehat{G}$,
\end{itemize}
then $|N|\leqslant |E|$.}\vspace{1mm}

\proof Consider the group ring element
$$X=E^2E^{(-1)}-\frac{|E|^2}{4}E\in\Z[G].$$ Since $\chi(X)=0$
for every non-principal character $\chi\in\widehat{G}$, $X=\lambda
G$ for some integer $\lambda$. Applying the principal character of
$G$, one gets $\lambda=(3|E|^3)/(4|G|)=(3|E|^2)/(4|N|)$ and $|N|$
divides $|E|^2/4$ as $G$ is a $2$-group. When $|E|=2$ or $4$, we
clearly have $|N|\leqslant |E|$. When $|E|=8$, if $\lambda$ is even,
then $|N|$ divides $|E|^2/8$ and therefore $|N|\leqslant |E|$. If
$\lambda$ is odd, then $\lambda=3$ and
$$E^2E^{(-1)}=\frac{|E|^2}{4}E+3G.$$
>From $E^2=E^{(2)}+2U$ for some $U\in\Z[G]$, one gets
$E^{(2)}E^{(-1)}-G=2V$ for some $V\in\Z[G]$. This implies that for
every $g\in G$, the coefficient of $g$ in $E^{(2)}E^{(-1)}$ is odd,
and therefore $|G|\leqslant |E^{(2)}E^{(-1)}|=|E|^2$ and
$|N|\leqslant |E|$.\qed

{\cor\label{II} Let $E$ be of type $\rm{(II)}$. If $E$ has the
property that $|\chi(E)|=|E|/2$ or $0$ for every non-principal
character $\chi\in\widehat{G}$, then $|N|\leqslant 8$.}\vspace{1mm}

\proof Since $E$ is a transversal of $N$ in $G$ and also it is the
union of $8$ cosets of the subgroup $H$, one must have $|H\cap N|=1$
and $E{'}$ is a transversal of $NH$ in $G$. Let $\overline{E}$ be
the image of $E$ in $\overline{G}=G/H$. Then
$\overline{E}=|H|\overline{E}{'}$, where $\overline{E}{'}$ is the
image of $E{'}$ in $\overline{G}$. Since $E{'}$ is a transversal of
$NH$ in $G$, its image $\overline{E}{'}$ contains $8$ elements and
is a transversal of $\overline{N}\cong N$ in $\overline{G}$. Because
$E$ has the property that $|\chi(E)|=|E|/2$ or $0$ for every
non-principal character $\chi\in\widehat{G}$, the transversal
$\overline{E}{'}$ also has the property that
$|\chi(\overline{E}{'})|=|\overline{E}{'}|/2$ or $0$ for every
non-principal character $\chi\in\widehat{\overline{G}}$. By
Lemma~\ref{size}, $|N|\leqslant 8$. \qed

If the set $E$ in Lemma~\ref{E2} is a transversal of a subgroup $N$
of $G$, then by a translation, we can assume that such a transversal
contains a subgroup $A$ of $G$ with $|A|=|E|/4$. Obviously $|A\cap
N|=1$ as $A$ is contained in a transversal of $N$, and $[G:NA]=4$.
The following two lemmas show that $E$ is of type (I) or (II) when
$G/NA\cong(\Z/2\Z)^2$ and is of type (I) when $G/NA\cong\Z/4\Z$ and
$N$ is elementary.

{\lem\label{2x2} Let $G$ be an abelian $2$-group of exponent no
larger than $4$ and $N$ a subgroup of $G$ such that $|N|\geqslant
8$. Let $E$ be a transversal of $N$ in $G$ which contains a subgroup
$A$ of $G$ such that $|A|=|E|/4$. If
\begin{itemize}
\item[(1)] for every character $\chi\in\widehat{G}$, the absolute value
$|\chi(E)|=|E|/2$ or $0$,  and
\item[(2)] $G/NA\cong(\Z/2\Z)^2$,
\end{itemize}
then $E$ must be of type $\rm{(I)}$ or $\rm{(II)}$.}\vspace{1mm}

\proof Let $\{1,x,y,xy\}$ be a transversal of $NA$ in $G$. Since
$G/NA\cong(\Z/2\Z)^2$, we have three index $2$ subgroups
$H_1=(1+x)NA$, $H_2=(1+y)NA$ and $H_3=(1+xy)NA$ of $G$. By
Lemma~\ref{cong}, the set $E\cap H_i$ has the property that
$\chi(E\cap H_i)\equiv 0$\\$(\!\!\!\!\mod |E|/4)$ for all $i$. Note
that $E$ is a transversal of $N$ in $G$ and $N<H_i$ for all $i$,
therefore $E\cap H_i$ is a transversal of $N$ in $H_i$ and $|E\cap
H_i|=|E|/2$ for all $i$. Since $E\cap NA=A$, one gets
$$|E\cap xNA|=|E\cap yNA|=|E\cap xyNA|=|E|/4$$ and $$\chi(E\cap
xNA)\equiv\chi(E\cap yNA)\equiv\chi(E\cap xyNA)\equiv
0\;\;(\!\!\!\!\!\!\mod |E|/4).$$ This implies that $E\cap xNA$,
$E\cap yNA$ and $E\cap xyNA$ are cosets of subgroups of order
$|E|/4$ in $G$. Hence we can assume that $E\cap xNA=xuA_1$, $E\cap
yNA=yvA_2$ and $E\cap xyNA=xywA_3$, where $A_i$ is a subgroup of
$NA$ such that $|A_i|=|E|/4$ and $|A_i\cap N|=1$ for $i=1$, $2$ and
$3$, and $u,v$ and $w$ are in $N$ as $NA_i=NA$ for all $i$, and
$E=A+xuA_1+yvA_2+xywA_3$.

If two of the four groups $A$, $A_1$, $A_2$ and $A_3$ are the same,
we claim that $E$ is of type (I). Without loss of generality, we
assume that $A=A_2$. Then $A_1=A_3$ because $\chi(E)\equiv 0$(mod
$|E|/2$) implies that $\chi$ is principal on even number of the
groups $A$, $A_1$, $A_2$ and $A_3$ for every $\chi\in\widehat{NA}$.
Therefor $A+yuA_2=A(1+yu)$ and $xvA_1+xywA_3=xvA_1(1+ywv^{-1})$. We
now want to show that $(yu)^2\in A$ or $(ywv^{-1})^2\in A_1$, and it
follows that $A(1+yu)$ is a group or $A_1(1+ywv^{-1})$ is a group,
and subsequently it implies that both of them are groups as
$\chi(E)\equiv 0$(mod $|E|/2$) for all $\chi\in\widehat{G}$, and
hence $E$ is of type (I). Suppose to the contrary. Then
$(yu)^2\notin A$ and $(ywv^{-1})^2\notin A_1$. Let's look at
$(yu)^2(ywv^{-1})^2$. Since $(yu)^2(ywv^{-1})^2=(uwv^{-1})^2\in N$
and $A\cap N=\{1\}$, if $(yu)^2(ywv^{-1})^2\ne 1$, then we can find
a $\chi\in\widehat{G}$ such that
$\chi((yu)^2(ywv^{-1})^2)=\chi((yu)^2)=-1$ and $\chi$ is principal
on $A$, which implies that
$$\chi(A(1+yu))=(1\pm\sqrt{-1})|E|/4\mbox{ and }
\chi(xvA_1(1+ywv^{-1}))\equiv 0\;\;(\!\!\!\!\mod |E|/2)$$ as
$\chi((ywv^{-1})^2)=1$. This contradicts the condition that
$\chi(E)\equiv 0\;\;(\!\!\!\!\mod |E|/2)$. So $(yu)^2(ywv^{-1})^2=
1$ and $(yu)^2=(ywv^{-1})^2$. Now $(yu)^2=(ywv^{-1})^2$ can not be
in $\<a,A_1\>$ because if $(yu)^2=(ywv^{-1})^2\in\<a,A_1\>$, then
there is a $\chi\in\widehat{G}$ such that $\chi((yu)^2)=-1$ and
$\chi$ is principal on $A$ as $(yu)^2\notin A$, which would ensure
$\chi(A_1)=0$ as $(yu)^2=(ywv^{-1})^2\in\<a,A_1\>$ and
$\chi(E)=(1\pm\sqrt{-1})|E|/4$, contradicting $\chi(E)\equiv 0$ (mod
$|E|/2$). Now since $(yu)^2=(ywv^{-1})^2\notin\<a,A_1\>$, one must
have $A=A_1$ because, again, if $A\ne A_1$, then there is a
$\chi\in\widehat{G}$ such that $\chi((yu)^2)=-1$, $\chi$ is
principal on $A$ but not principal on $A_1$ as $A\ne A_1$, which
implies that $\chi(E)=(1\pm\sqrt{-1})|E|/4$ contradicting
$\chi(E)\equiv 0$ (mod $|E|/2$). Hence $E=A(1+xu+yv+xyw)$. In the
group $\overline{G}=G/A$, the image of $\{1,xu,yv,xyw\}$ form a
transversal of $\overline{N}\cong N$ in $\overline{G}$ and has the
property that
$|\chi(1+\overline{xu}+\overline{yv}+\overline{xyw})|=2$ or $0$ for
all non-principal character $\chi\in\widehat{\overline{G}}$. By
Lemma~\ref{size}, $N\leqslant 4$, contradicting $N\geqslant 8$.

If no two of $A$, $A_1$, $A_2$ and $A_3$ are the same, let $A_0=A$,
we claim that
$$|\bigcap_{i=0}^{3}A_{i}|=|E|/8.$$  Let $H=\bigcap_{i=0}^{3}A_{i}.$
By taking the quotient $G/H$, we may assume that $H=\{1\}$ and we
need to show that $|E|=8$. From $|H|=1$, we can see that the
intersections of any three of the four groups $A$, $A_1$, $A_2$ and
$A_3$ are trivial.  This is because if there are three groups, say
$A_1$, $A_2$ and $A_3$ such that $\bigcap_{i=1}^{3}A_i\ne\{1\}$,
then there is a character $\chi\in\widehat{G}$ such that $\chi$ is
non-principal on $\bigcap_{i=1}^{3}A_i$ and principal on $A_0=A$ as
$A\cap(\bigcap_{i=1}^{3}A_i)=\{1\}$, which implies
$|\chi(E)|=|E|/4$, a contradiction. We now calculate the sizes of
pairwise intersections of the $4$ groups. Let
$F=NA_{0}=NA_{1}=NA_{2}=NA_{3}$ and consider the images
$\overline{A}_{1}$, $\overline{A}_{2}$ and $\overline{A}_{3}$ of
$A_{1}$, $A_{2}$ and $A_{3}$ respectively in the group $F/A_{0}\cong
N$. None of these three images is trivial as $A_{i}\nleqslant A_{0}$
for all $i\ne 0$. Also every character $\chi\in\widehat{F/A_{0}}$
must be trivial on at least one of the three images as such a $\chi$
of $F/A_{0}$ can be extended to a character of $G$ which is
principal on $A_{0}$. If there are two images, say $\overline{A}_1$
and $\overline{A}_2$, such that
$\overline{A}_1\cap\overline{A}_2\ne\{1\}$, then for any $1\ne
u\in\overline{A}_1\cap\overline{A}_2$ and $1\ne v\in\overline{A}_3$,
by Lemma~\ref{2elts}, there is a character
$\chi\in\widehat{F/A_{0}}$ such that $\chi(u)\ne 1$ and $\chi(v)\ne
1$, which implies there are characters which are non-principal on
all three images and this is not possible. Hence the pairwise
intersections of the three images must be trivial. Let $1\ne
a_i\in\overline{A_{i}}$ for $i=1,2$ and $3$. Since there is no
character $\chi$ such that $\chi(a_i)\ne 1$ for all $i$, the set
$\{1, a_1,a_2,a_3\}$ must be a group isomorphic to $(\Z/2\Z)^2$. Fix
$a_1$ and $a_2$, we can see that $\overline{A}_{3}=\{1,a_3\}$.
Similarly, $\overline{A}_{i}=\{1,a_i\}$ for all $i=1$, $2$ and $3$.
Hence $|A_{i}/(A_{0}\cap A_{i})|=|\overline{A}_{i}|=2$ for all $i$,
and similarly, $|A_{i}/(A_{i}\cap A_{j})|=2$ for all $i\ne j$. Since
$|(A_{i}\cap A_{j})\cap(A_{i}\cap A_{k})|=1$ when $i$, $j$ and $k$
are distinct, either $A_{i}\cong(\Z/2\Z)^2$ for all $i$, or
$A_{i}\cong\Z/2\Z$ for all $i$ as these are the only groups
containing index $2$ subgroups with trivial intersections. If
$A_{i}\cong(\Z/2\Z)^2$, then
$\<a_{0},A_{1},A_{2},A_{3}\>\cong(\Z/2\Z)^3$, which is impossible as
there are characters that are trivial on exactly one of these
$A_i$'s since they are hyperplanes of $\<a_{0},A_{1},A_{2},A_{3}\>$.
So $|A_{i}|=2$ for $i=0$, $1$, $2$ and $3$ and $|E|=8$.\qed

{\lem\label{Z4} Let $G$ be an abelian $2$-group of exponent no
larger than $4$ and $N$ an elementary abelian subgroup of $G$ such
that $|N|\geqslant 8$. Let $E$ be a transversal of $N$ in $G$ which
contains a subgroup $A$ of $G$ such that $|A|=|E|/4$. If
\begin{itemize}
\item[(1)] for every character $\chi\in\widehat{G}$, the absolute value
$|\chi(E)|=|E|/2$ or $0$,  and
\item[(2)] $G/NA\cong\Z/4\Z$,
\end{itemize}
then $E$ must be of type $\rm{(I)}$.}\vspace{1mm}

\proof If $G/NA\cong\Z/4\Z$, we may further assume that
$E=A+xA_1+x^2A_2+x^3A_3=(A+x^2A_2)+x(A_1+x^2A_3)$, where $A_1$,
$A_2$ and $A_3$ are transversals of $N$ in $NA$ and $x$ is an
element of order $4$ in $G$ such that $|\<x\>\cap (NA)|=1$. Since
$A+x^2A_2=E\cap\<x^2\>NA$, By Lemma~\ref{cong},
$\chi(A+x^2A_2)\equiv 0$ $(\!\!\!\!\mod |E|/4)$ for all
$\chi\in\widehat{\<x^2\>NA}$. It follows that $A_2=uB$ for some
subgroup $B$ of order $|E|/4$ in $NA$ with $|B\cap N|=1$ and
$NA=NB$, and $u\in N$ because $A$ is a group of order $|E|/4$. We
now show that $A=B$. Let $F=NA=NB$. If $A\ne B$, then there is a
$\chi\in\widehat{G}$ such that $\chi(A)\ne\chi(B)$ and
$\chi(x^2)=-1$ as $|\<x\>\cap F|=1$. Since $F/A\cong F/B\cong N$ is
elementary abelian and either $A$ or $B$ is contained in the kernel
of $\chi$, the character $\chi|_F$ must be of order $2$. Therefore
$\chi(A_1)$ and $\chi(A_3)$ are real, and consequently,
$\chi(A_1+x^2A_3)$ is real, which contradicts $|\chi(E)|=|E|/2$ or
$0$ for all non-principal $\chi\in\widehat{G}$. Hence $A=B$ and
$A+x^2A_2=A(1+ux^2)$ is a group as $ux^2$ is of order $2$. Now for
every $\chi\in\widehat{G}$,
$$\chi(x(A_1+x^2A_3))=\chi(E)-\chi(A+x^2A_2)\equiv
0\;\;(\!\!\!\!\!\!\mod |E|/2).$$ This implies that $A_1+x^2A_3$ is a
coset of a subgroup of order $|E|/2$ in $G$ and $E$ is of type (I).
\qed

We now change the conditions imposed on $G/NA$ to conditions on
characters and obtained following two corollaries which are actually
the ones used in the proof of the main theorem.

{\cor\label{EI} Let $G$ be an abelian $2$-group of exponent no
larger than $4$ and $N$ an elementary abelian subgroup of $G$ such
that $|N|\geqslant 8$. If a transversal $E$ of $N$ in $G$ satisfies
\begin{itemize}
\item[(1)] for every non-principal character $\chi\in\widehat{G}$, the
absolute value $|\chi(E)|=|E|/2$ or $0$,  and
\item[(2)] there is a character $\phi\in\widehat{G}$ of order $2$ such
that $|\phi(E)|=|E|/2$,
\end{itemize}
then $E$ must be of type $\rm{(I)}$ or $\rm{(II)}$.}\vspace{1mm}

\proof By Lemma~\ref{E2}, we may assume that the set $E$ contains a
subgroup $A$ of order $|E|/4$ in $G$. Since $E$ is a transversal of
$N$ in $G$, $|A\cap N|=1$ and $[G:NA]=4$. By Lemma~\ref{2x2} and
\ref{Z4}, the transversal $E$ is either of type $\rm{(I)}$ or
$\rm{(II)}$.\qed

{\cor\label{EJ} Let $G$ be an abelian $2$-group of exponent $4$ and
$N$ a subgroup of $G$ such that $|N|\geqslant 8$ and the exponent of
$N$ is $4$. If a transversal $E$ of $N$ in $G$ satisfies
\begin{itemize}
\item[(1)] for every non-principal character $\chi\in\widehat{G}$,
the absolute value $|\chi(E)|=|E|/2$ or $0$,  and
\item[(2)] there is a character $\phi\in\Phi(\widehat{G})$ such
that $|\phi(E)|=|E|/2$,
\end{itemize}
then $E$ must be of type $\rm{(I)}$ or $\rm{(II)}$.}\vspace{1mm}

\proof Let $K={\rm Ker}(\phi)$ be the kernel of $\phi$. Since
$\phi\in\Phi(\widehat{G})$, $\phi$ is of order $2$ and $|G/K|=2$. By
the proof of Lemma~\ref{E2}, we can assume that $E\cap K$ is a
subgroup of $K$ of order $|E|/4$. We know that $G=K+xK$ for some
$x\notin K$. the element $x$ is of order $4$ as
$\phi\in\Phi(\widehat{G})$ and we can also assume that $x\in N$ as
$E$ is a transversal of $N$ in $G$. We now show that there is an
index $2$ subgroup $H$ of $K$ such that $|\<x^2\>\cap H|=1$. This is
equivalent to show that $x^2\notin\Phi(K)$, i.e. there is no $y\in
K$ such that $y^2=x^2$. Suppose that there is an $y\in K$ such that
$y^2=x^2$. Then $(xy)^2=1$. This implies that $xy\in K$ as
$\phi\in\Phi(\widehat{G})$, and hence $x\in K$, a contradiction. Let
$H$ be a subgroup of $K$ such that $|\<x^2\>\cap H|=1$ and
$K=\<x^2\>H$. The group $E\cap K$ is not contained in $H$ because if
$E\cap K\leqslant H$, then there is a character $\chi\in\widehat{G}$
such that $\chi$ is principal on $H$ and $\chi(x^2)=-1$, which
implies that
$$\chi(E)=\chi(E\cap K)+\chi(x)\chi(x^{-1}E\cap
K)=(|E|/4)\pm\chi(x^{-1}E\cap K)\sqrt{-1}$$ and $\chi(x^{-1}E\cap
K)$ is real as $\chi|_K$ is of order $2$ since $H\leqslant{\rm
Ker}(\chi)$, contradicting the condition that $|\chi(E)|=|E|/2$ or
$0$ for every non-principal character $\chi\in\widehat{G}$. Hence
$E\cap H$ is an index $2$ subgroup of $E\cap K$ and $E\cap K=(E\cap
H)(1+x^2y)$ for some $y\in H$ such that $y\notin E$ but $y^2\in E$.
Let $F=(E\cap H)(1+y)$. $F$ is a subgroup of $H$ of order $|E|/4$.
Let $H_1$ be an index $2$ subgroup of $H$ such that $E\cap
H\leqslant H_1$. If $y\notin H_1$, then there is a
$\chi\in\widehat{G}$ such that $\chi$ is principal on $H_1$ but not
on $H$ and $\chi(x^2)=-1$. Such a $\chi$ is real on $K$ as
$K=\<x^2\>H$ and $\chi|_H$ is real. Hence, again
$$\chi(E)=\chi(E\cap K)+\chi(x)\chi(x^{-1}E\cap
K)=(|E|/4)\pm\chi(x^{-1}E\cap K)\sqrt{-1}$$contradicting the
condition that $|\chi(E)|=|E|/2$ or $0$ for every non-principal
character $\chi\in\widehat{G}$. Therefore $y\in H_1$ for every index
$2$ subgroup $H_1$ of $H$ that contains $E\cap H$. This is
equivalent to that the rank of $H/(E\cap H)$ is equal to that of
$H/(E\cap H)(1+y)$, or the image of $y$ in $H/(E\cap H)$ is in
$\Phi(H/(E\cap H))$, i.e. $y=st^2$ for some $s\in E\cap H$ and $t\in
H$ but $|(E\cap H)\cap\<t\>|=1$. Therefore $F=(E\cap H)(1+y)=(E\cap
H)(1+t^2)$ and $E\cap K=(E\cap H)(1+x^2t^2)$. Now $t^2\notin N(E\cap
H)$ because if $t^2=uv$ with $u\in N$ and $v\in E\cap H$, then
$$E\cap K=(E\cap H)(1+x^2t^2)=(E\cap H)(1+x^2uv)=(E\cap
H)(1+x^2u),$$ which is impossible as $E$ is a transversal of $N$ in
$G$ and $x^2u\in N$. Hence $|\<t\>\cap N(E\cap H)|=1$ and $G/N(E\cap
H)\cong\Z/2\Z\times\Z/4\Z$. Since $N(E\cap K)=N(E\cap
H)(1+x^2t^2)=N(E\cap H)(1+t^2)$ as $x^2\in N$, one gets $G/N(E\cap
K)\cong(\Z/2\Z)^2$ and, by Lemma~\ref{2x2}, $E$ is of type (I) or
(II).\qed

We finish this section with a lemma that shows the groups involved
in the type (I) transversals of $(2^{m+1},2^{m},2^{m-1})$-building
sets form some kind of orthogonal system.

{\lem\label{EIJ} Let $G$ be an abelian group of order $2^{2m+1}$ and
of exponent at most $4$ and $N$ a subgroup of $G$ of order $2^m$.
Let $E_1=a_1H_1+a_2H_2$, $E_2=a_3H_3+a_4H_4$, $\cdots$,
$E_k=a_{2k-1}H_{2k-1}+a_{2k}H_{2k}$ be a family of type $\rm{(I)}$
transversals of $N$ in $G$, where $H_i$ are subgroups of order $2^m$
in $G$ and $H_i\cap N=\{1\}$ for all $i=1$, $2$, $\cdots$, $k$. If
for every non-principal character $\chi\in\widehat{G}$, there is at
most one $E_i$ such that $|\chi(E_i)|=2^m$, and $\chi(E_j)=0$ for
all $j\ne i$, then $H_s\cap H_t=\{1\}$ and $NH_s=NH_t=H_sH_t$ for
all $1\leqslant s\ne t\leqslant 2k$.}\vspace{1mm}

\proof According to Lemma~\ref{I} and the condition that for every
non-principal character $\chi\in\widehat{G}$, there is at most one
$E_i$ such that $|\chi(E_i)|=2^m$, and $\chi(E_j)=0$ for all $j\ne
i$, we can see that $\chi$ can be principal on at most one $H_i$,
$1\leqslant i\leqslant k$ when $\chi$ is non-principal on $N$. If
$\chi$ is non-principal on $N$ and principal on $H_s$ for some
$1\leqslant s\leqslant k$, then $\chi$ is a character of
$\overline{G}=G/H_s$. The image $\overline{N}$ of $N$ in
$\overline{G}$ is isomorphic to $N$ as $H_s\cap N=\{1\}$. For every
$t\ne s$, the group $\overline{N}$ must be contained in the image
group $\overline{H}_t$ of $H_t$ because otherwise there is a
$\chi\in\widehat{\overline{G}}$ such that $\chi$ is principal on
$\overline{H}_t$ and non-principal on $\overline{N}$. Therefore
$\overline{H}_t=\overline{N}$ and $H_tH_s=NH_s$ for all $1\leqslant
s\ne t\leqslant 2k$.\qed

\section{Proof of Theorem~\ref{main}}

In this we prove Theorem~\ref{main}. We start with two simple lemmas
and their proofs, which we omit, are easy exercises of applications
of characters.

{\lem\label{quotient} If an abelian group $G$ of order $2^{m+1}n$
contains $(2^{m+1},2^m,2^{m-1})$-building sets relative to a
subgroup $N$ of order $n$ in $G$, then for any subgroup $H$ of $N$,
the group $G/H$ contains $(2^{m+1},2^m,2^{m-1})$-building sets
relative to a subgroup $N/H$.}\vspace{1mm}

{\lem\label{number} Let $G$ be a group of order $2^{m+n+1}$ and $N$
a subgroup of $G$ of order $2^n$. If $E_1$, $E_2$, $E_3$, $\cdots$,
$E_{2^{m-1}}$ form $(2^{m+1},2^m,2^{m-1})$-building sets in $G$
relative to $N$, then
$$\sum_{k=1}^{2^{m-1}}E_kE_k^{(-1)}=2^{2m}+2^{2m-n}(G-N).$$
}\vspace{1mm}

The next lemma excludes type (II) transversals from any
$(2^{m+1},2^m,2^{m-1})$-building sets in a group of order
$2^{2m+1}$.

{\lem\label{noII} Let $G$ be a group of order $2^{2m+1}$ and
$m\geqslant 3$. There are no type ${\rm (II)}$ transversals in any
$(2^{m+1},2^m,2^{m-1})$-building sets relative to any subgroup of
order $2^m$ in $G$.}\vspace{1mm}

\proof If $E_1$, $E_2$, $\cdots$, $E_{2^{m-1}}$ form
$(2^{m+1},2^m,2^{m-1})$-building sets relative to a subgroup of
order $2^m$ in $G$ and, say, $E_1$ is of type (II), then $E_1=HE{'}$
for some non-trivial subgroup $H$ of $G$ and some subset $E{'}$ of
$G$ with $|E{'}|=8$. Let $x\in H$ and $x\ne 1$. Then the coefficient
of $x$ in $E_1E_1^{(-1)}$ is $|E_1|=2^{m+1}$, which contradicts
Lemma~\ref{number}.\qed

We now prove Theorem~\ref{main}.  \vspace{2mm}

\noindent{\it Proof of Theorem}~\ref{main}. The existence part of
the theorem is given in \cite{Dillon0, Mc}. By Lemma~\ref{bset} and
results in \cite{MS1,MS2}, we need to prove that any abelian
$2$-group $G$ of order $2^{2m+1}$ and of exponent no larger than $4$
containing $(2^{m+1},2^m,2^{m-1})$-building sets relative to a
subgroup $N$ of order $2^m$ must contain an elementary abelian
subgroup of index $2$ in $G$.  Let $E_1$, $E_2$, $\cdots$,
$E_{2^{m-1}}$ be $(2^{m+1},2^m,2^{m-1})$-building sets relative to
$N$ in $G$. We are going to show that $N$ is an elementary abelian
subgroup of $G$.

Suppose $N$ is not elementary abelian. Then there are characters
$\chi\in\Phi(\widehat{G})$ which are not principal on $N$. By the
definition of building sets, $E_1$, $E_2$, $\cdots$, $E_{2^{m-1}}$
are transversals of $N$ in $G$. By Corollary~\ref{EJ} and
Lemma~\ref{noII}, there are type (I) transversals among $E_1$,
$E_2$, $\cdots$, $E_{2^{m-1}}$. Assume that $E_1=a_1H_1+a_2H_2$,
$E_2=a_3H_3+a_4H_4$, $\cdots$, $E_k=a_{2k-1}H_{2k-1}+a_{2k}H_{2k}$,
$k\leqslant 2^{m-1}$ are of type (I) transversals and the rest are
not.
\vspace{1.0mm}
\newline{\bf Claim 1:} $k=|\Phi(N)|/2$.
\vspace{1.0mm}
\newline Let $H=NH_1$. By Lemma~\ref{EIJ},
$H=NH_s=H_sH_t$ for all $1\leqslant s\ne t\leqslant 2k$ and $H$ is
an index $2$ subgroup of $G$. Let $H_0=N$. Since $H_i<H$ and $H\cong
H_i\times H_i$ for all $0\leqslant i\leqslant 2k$,
$\Phi(H_i)<\Phi(H)$ and $\Phi(H)\cong\Phi(H_i)\times\Phi(H_i)$ for
all $0\leqslant i\leqslant 2k$. Since $\Phi(H_i)\cong\Phi(N)$ as
$H_i\cong N$ and
$(\Phi(H_i)\setminus\{1\})\cap(\Phi(H_j)\setminus\{1\})=\emptyset$
as $H_i\cap H_j=\{1\}$ for all $0\leqslant i\ne j\leqslant 2k$,
there can be at most
$|\Phi(H)\setminus\{1\}|/|\Phi(N)\setminus\{1\}|$ such subgroups of
$H_i$, i.e. $2k+1\leqslant
(|\Phi(H)|-1)/(|\Phi(N)|-1)=(|\Phi(N)|^2-1)/(|\Phi(N)|-1)=|\Phi(N)|+1.$
Therefore $k\leqslant|\Phi(N)|/2$. On the other hand, for each type
(I) transversal $E_i$, a $\chi\in\Phi(\widehat{G})$ has
$|\chi(E_i)|=2^m$ if and only if $\chi$ is non-principal on $N$,
principal on the maximum elementary abelian subgroup of $G$ and
principal on exactly one of $H_{2i-1}$ and $H_{2i}$. There are
exactly $2[(|\Phi(G)|/|\Phi(N)|)-(|\Phi(G)|/|\Phi(N)|^2)]$ such
$\chi\in\Phi(\widehat{G})$. The group $\Phi(\widehat{G})$ contains
$|\Phi(G)|-(|\Phi(G)|/|\Phi(N)|)$ characters which are not principal
on $N$. Therefore, by By Corollary~\ref{EJ} and Lemma~\ref{noII},
there should be at least
$$\frac{|\Phi(G)|-(|\Phi(G)|/|\Phi(N)|)}{2[(|\Phi(G)|/|\Phi(N)|)-(|\Phi
(G)|/|\Phi(N)|^2)]}=
\frac{|\Phi(N)|}{2}$$ $E_i$'s of type (I). Hence $k=|\Phi(N)|/2$.
\vspace{1.0mm}
\newline{\bf Claim 2: $k\leqslant 2$.}
\vspace{1.0mm}
\newline Since $k=|\Phi(N)|/2$, the sets $\Phi(H_i)\setminus\{1\}$,
$i=0$, $1$,
$\cdots$, $|\Phi(N)|$ form a partition of $\Phi(H)\setminus\{1\}$.
By Lemma~\ref{Eg}, for each $|\Phi(N)|/2\leqslant i\leqslant |N|/2$,
$E_i$ contains a coset of a subgroup, say $K_i$, of order $2^{m-1}$
of $G$. By Lemma~\ref{number}, $K_i\cap\Phi(H)=\{1\}$. Therefore
$K_i\cap H$ is an elementary abelian group. Since $K_i\cap N=\{1\}$
as $E_i$ is a transversal of $N$ in $G$, one has
$$\Rank(H)\geqslant\Rank(K_i\cap H)+\Rank(N)=\Rank(K_i\cap
H)+[\Rank(H)/2].$$ Therefore $\Rank(H)\geqslant 2\Rank(K_i\cap H)$.
Since $H$ is of index $2$ in $G$ and $|K_i|=2^{m-1}$, the order
$|K_i\cap H|\geqslant 2^{m-2}$ and the rank $\Rank(K_i\cap
H)\geqslant m-2$. This implies $\Rank(H)\geqslant 2(m-2)=2m-4$ and
$\Rank(\Phi(H))\leqslant 4$ as $|H|=2^{2m}$. Hence
$k=|\Phi(N)|/2=\sqrt{|\Phi(H)|}/2\leqslant 2$.
\vspace{1.0mm}
\newline{\bf Claim 3: $k=1$.}
\vspace{1.0mm}
\newline If $k=2$, then
$N\cong(\Z/4\Z)^2\times(\Z/2\Z)^{m-4}$,
$H\cong(\Z/4\Z)^4\times(\Z/2\Z)^{2m-8}$. Also for each $i>2$, we
have $K_i\cap H\cong(\Z/2\Z)^{m-2}$ and $N(K_i\cap
H)\cong(\Z/4\Z)^2\times(\Z/2\Z)^{2m-6}$.  Let $x\in K_i$ such that
$x\notin H$. Then
$$G/NK_i=\<H,x\>/\<N(K_i\cap H),x\>\cong H/N(K_i\cap H)\cong(\Z/2\Z)
^2.$$ By
Lemma~\ref{2x2} and \ref{noII}, $E_i$ is of type (I), contradicting
the assumption that $E_i$ is not of type (I) for
$|\Phi(N)|/2\leqslant i\leqslant |N|/2$. Hence $k=1$.

Since $k=|\Phi(N)|/2=1$, one has
$$N\cong\Z/4\Z\times(\Z/2\Z)^{m-2},\mbox{ and }
H\cong(\Z/4\Z)^2\times(\Z/2\Z)^{2m-4}.$$  Also by Lemma~\ref{2x2}
and \ref{noII}, $G/NK_i\cong\Z/4\Z$ for $i\geqslant 2$.  From
$|G/H|=2$, one has either $\Rank(\Phi(G))=3$ and
$G\cong(\Z/4\Z)^3\times(\Z/2\Z)^{2m-5}$, or $\Rank(\Phi(G))=2$ and
$G\cong(\Z/4\Z)^2\times(\Z/2\Z)^{2m-3}$.
\vspace{1.0mm}
\newline{\bf Claim 4: $\Rank(\Phi(G))\ne 3$.}
\vspace{1.0mm}
\newline If $\Rank(\Phi(G))=3$ and
$G\cong(\Z/4\Z)^3\times(\Z/2\Z)^{2m-5}$, let's consider the group
$\overline{G}=G/\Phi(N)\cong(\Z/4\Z)^2\times(\Z/2\Z)^{2m-4}$. In
$\overline{G}$, the image $\overline{N}=N/\Phi(N)$ of $N$ is an
elementary abelian group of order $2^{m-1}$, and, by
Lemma~\ref{quotient}, the images $\overline{E}_1$, $\overline{E}_2$,
$\cdots$, $\overline{E}_{2^{m-1}}$ of $E_1$, $E_2$, $\cdots$,
$E_{2^{m-1}}$ form $(2^{m+1},2^m,2^{m-1})$-building sets relative to
$\overline{N}$ in $\overline{G}$. The transversal $\overline{E}_1$
is of type (I) since $E_1$ is. For $i\geqslant 2$, $\overline{E}_i$
contains the group $\overline{K}_i$, the image of $K_i$ in $G$, and
$\overline{K}_i\cong K_i$ as $K_i\cap N=\{1\}$. Since
$\overline{G}/\overline{N}\overline{K}_i\cong G/NK_i\cong\Z/4\Z$ and
$\overline{N}$ is elementary abelian, by Lemma~\ref{Z4}, the
transversal $\overline{E}_i$ is of type (I) for every $i\geqslant
2$. We can now assume that
$\overline{E}_i=a{'}_{2i-1}H{'}_{2i-1}+a{'}_{2i}H{'}_{2i}$ for some
subgroups $H{'}_{2i-1}$ and $H{'}_{2i}$ of order $2^m$ in
$\overline{G}$ and some $a{'}_{2i-1}$ and $a{'}_{2i}$ in
$\overline{G}$ for $i=1$, $2$, $\cdots$, $2^{m-1}$. The group
$H{'}_j\cap\overline{N}=\{1\}$ for all $j$ as $\overline{E}_i$ is a
transversal of $\overline{N}$. Hence $H{'}_j$ is not elementary
abelian, i.e. $\Phi(H{'}_j)$ is non-trivial, for all $j$ because
$\Rank(\overline{G})=2m-2$ while $|H{'}_j|=2^m$,
$|\overline{N}|=2^{m-1}$ and $\overline{N}$ is elementary abelian.
Also by Lemma~\ref{I}, $H{'}_{2i-1}\cap H{'}_{2i}$ is of order $2$
and contains $\Phi(H{'}_{2i-1})$ and $\Phi(H{'}_{2i})$, and
therefore $H{'}_{2i-1}\cap H{'}_{2i}=\Phi(H{'}_{2i-1})=
\Phi(H{'}_{2i})$ is contained in $\Phi(\overline{G})$ for all $i$.
By Lemma~\ref{number}, $(H{'}_{2s-1}\cap
H{'}_{2s})\ne(H{'}_{2t-1}\cap H{'}_{2t})$ for all $s\ne t$, i.e.
$H{'}_{2i-1}\cap H{'}_{2i}$ are distinct subgroups of order $2$ in
$\Phi(\overline{G})$. Since $|\Phi(\overline{G})|=4$, the group
$\Phi(\overline{G})$ only has three subgroups of order $2$. Hence
$2^{m-1}\leqslant 3$ and $m\leqslant 2$, contradicting $m\geqslant
3$.
\vspace{1.0mm}
\newline{\bf Claim 4: $\Rank(\Phi(G))\ne 2$.}
\vspace{1.0mm}
\newline If $\Rank(\Phi(G))=2$ and $G\cong(\Z/4\Z)^2\times(\Z/2\Z)^{2m-
3}$,
let $N{'}$ be an elementary abelian subgroup of $N$ of order
$2^{m-2}$ such that $N{'}\cap\Phi(N)=\{1\}$. By
Lemma~\ref{quotient},  the images $\overline{E}_1$,
$\overline{E}_2$, $\cdots$, $\overline{E}_{2^{m-1}}$ of $E_1$,
$E_2$, $\cdots$, $E_{2^{m-1}}$ in $\overline{G}=G/N{'}$ form
$(2^{m+1},2^m,2^{m-1})$-building sets relative to
$\overline{N}=N/N{'}\cong\Z/4\Z$. Moreover the transversal
$\overline{E}_1$ of $\overline{N}$ in $\overline{G}$ is of type (I)
since $E_1$ is and
$\overline{E}_1=\overline{a}_1\overline{H}_1+\overline{a}_2\overline{H}
_2$,
where $\overline{a}_1$, $\overline{a}_2$, $\overline{H}_1$ and
$\overline{H}_2$ are images of $a_1$, $a_2$, $H_1$ and $H_2$
respectively in $\overline{G}$. We now show that intersection
$\overline{H}_1\cap\overline{H}_2$ is an elementary abelian group of
order $2^{m-2}$. Since $H=H_1H_2=NH_1=NH_2$ and $N\cap H_1=N\cap
H_2=H_1\cap H_2=\{1\}$, there are isomorphisms $f_1:N\rightarrow
H_1$ and $f_2:N\rightarrow H_2$ such that $N=\{f_1(x)f_2(x)\in
H_1H_2=H~|~x\in N\}$. Hence
\begin{eqnarray*}
N{'}H_1&=&\{h_1f_2(y)\in H_1H_2=H~|~h_1\in H_1, y\in N{'}\},\\
N{'}H_2&=&\{f_1(x)h_2\in H_1H_2=H~|~x\in N{'}, h_2\in H_2\}
\end{eqnarray*} and
$N{'}H_1\cap N{'}H_2=\{f_1(x)f_2(y)\in H_1H_2=H~|~x,y\in N{'}\}\cong
N{'}\times N{'}$. This implies that
$$\overline{H}_1\cap\overline{H}_2=(N{'}H_1\cap N{'}H_2)/N{'}\cong
N{'}$$ is an elementary abelian group of order $2^{m-2}$. Let
$F=\overline{H}_1\cap\overline{H}_2$ and $\overline{N}=\<z\>$. By
Lemma~\ref{I}, we can assume that there is an element $x$ of order
$4$ in $\overline{G}$ such that $\overline{H}_1=\<F,x\>$ and
$\overline{H}_2=\<F,xz\>$. The group
$\overline{H}=\<\overline{H}_1,\overline{H}_2\>=\<F,x,z\>$ is an
index $2$ subgroup of $\overline{G}$. Let
$y\in\overline{G}\setminus\overline{H}$. Then $y^2=1$. Let
$K=\<F,x,y\>$. Then $K\cap\overline{N}=\{1\}$. For each $2\leqslant
i\leqslant 2^{m-1}$, the set $\overline{E}_i$ must have the form
$$\overline{E}_i=F_{i0}+F_{i1}z+F_{i2}z^2+F_{i3}z^3,$$ where
$F_{i0}$, $F_{i1}$, $F_{i2}$ and $F_{i3}$ are disjoint subsets of
$K$ such that $F_{i0}+F_{i1}+F_{i2}+F_{i3}=K$. Since for any
character $\chi\in\widehat{\overline{G}}$ with $\chi(z)=-1$,
$|\chi(\overline{E}_1)|=2^m$ if and only if $\chi$ is principal on
$\<F, x^2\>$, one finds that $F_{20}+F_{22}$, $F_{30}+F_{32}$,
$\cdots$, $F_{2^{m-1}0}+F_{2^{m-1}2}$ form $(2^m$, $2^{m-1}$,
$2^{m-1}-1)$-building sets in $K$ relative to $\<F,x^2\>$, and hence
we can further assume that
$$F_{i0}+F_{i2}=A_{i0}+A_{i2}x$$
for some subsets $A_{i0}$ and  $A_{i2}$ of sizes $2^{m-1}$ in
$\<y,F,x^2\>$ for all $i=2$, $3$,  $\cdots$, $2^{m-1}$. Clearly,
$F_{i1}+F_{i3}$ is the complement of $F_{i0}+F_{i2}$ in $K$,
therefore $F_{i1}+F_{i3}=A_{i1}+A_{i3}x$, where
$A_{i1}=\<y,F,x^2\>-A_{i0}$ and $A_{i3}=\<y,F,x^2\>-A_{i2}$ for all
$i\geqslant 2$. Let
$$\overline{E}_i=B_{i0}+B_{i1}z+B_{i2}x+B_{i3}xz,$$ where
\begin{eqnarray*} B_{i0}&=&\overline{E}_i\cap\<y,F,x^2,z^2\>,\\
B_{i1}&=&z^{-1}\overline{E}_i\cap\<y,F,x^2,z^2\>,\\
B_{i2}&=&x^{-1}\overline{E}_i\cap\<y,F,x^2,z^2\>,\\
B_{i3}&=&(xz)^{-1}\overline{E}_i\cap\<y,F,x^2,z^2\>.
\end{eqnarray*}
If $p:\<y,F,x^2,z^2\>\rightarrow\<y,F,x^2\>$ is the natural
projection, then $p(B_{is})=A_{is}$ for all $i\geqslant 2$ and
$s=0,1,2,3$, i.e. $B_{i0}$, $B_{i1}$, $B_{i2}$ and $B_{i3}$ are the
lift-backs of $A_{i0}$, $A_{i1}$, $A_{i2}$ and $A_{i3}$ from
$\<y,F,x^2\>$ to $\<y,F,x^2,z^2\>$. Clearly,
$B_{i0}+B_{i2}x=F_{i0}+F_{i2}z^2$ and
$B_{i1}+B_{i3}x=F_{i1}+F_{i3}z^2$ for all $i\geqslant 2$. Now we
want to evaluate the coefficients of $x^2z^2$ in
$\overline{E}_i\overline{E}_i^{(-1)}$ for $i\geqslant 2$. If the
coefficient of $x^2z^2$ in $\overline{E}_i\overline{E}_i^{(-1)}$ is
not $0$ for some $i\geqslant 2$, then $\overline{E}_i$ contains a
coset of $\<x^2z^2\>$. Obviously, such a coset must be in one of
$B_{is}$, for some $0\leqslant s\leqslant 3$. Without loss of
generality, we can assume that $B_{i0}$ contains a coset of
$\<x^2z^2\>$. Let $\chi\in\widehat{\<z^2\>}$ such that
$\chi(z^2)=-1$. Then $\chi$ induces a $\C$-algebra homomorphism,
which we denoted again by $\chi$,
$\chi:\C[K\times\<z^2\>]\rightarrow\C[K]$. Let
$\rho:\C[K]\rightarrow\C[K/\<x^2\>]$ be the $\C$-algebra
homomorphism induced by natural projection $K\rightarrow K/\<x^2\>$.
Then $\rho\chi(F_{i0}+F_{i2}z^2)\in\C[K/\<x^2\>]$. For every
character $\phi\in\widehat{K/\<x^2\>}$, $\phi\rho\chi$ is a
character of $K\times\<z^2\>$ with $\phi\rho\chi(z^2)=-1$ and
$\phi\rho\chi(x^2)=1$. It can be extended to $\overline{G}$ and
$$\phi\rho\chi(\overline{E}_i)=\phi\rho\chi(F_{i0}+F_{i2}z^2)\pm
\phi\rho\chi(F_{i1}+F_{i3}z^2)\sqrt{-1}.$$ Since
$\phi\rho\chi(F_{i0}+F_{i2}z^2)$ and
$\phi\rho\chi(F_{i1}+F_{i3}z^2)$ are real as $\phi\rho\chi(x^2)=1$,
one has $|\phi\rho\chi(F_{i0}+F_{i2}z^2)|=2^m$ or $0$. On the other
hand, $F_{i0}+F_{i2}z^2$ contains cosets of $\<x^2z^2\>$ and they
are in the kernel of $\rho\chi$. Therefore
$|\phi\rho\chi(F_{i0}+F_{i2}z^2)|\leqslant|F_{i0}+F_{i2}z^2|-2(\mbox
{number
of cosets of }\<x^2z^2\>)<2^m.$ Hence
$\phi(\rho\chi(F_{i0}+F_{i2}z^2))=0$ for all
$\phi\in\widehat{K/\<x^2\>}$. This implies that
$\rho\chi(F_{i0}+F_{i2}z^2)=0$ in $\C[K/\<x^2\>]$, or equivalently,
$\rho(F_{i0})=\rho(F_{i2})$. The fact that $F_{i0}$ and $F_{i2}$ are
disjoint subsets of $K$ implies that $F_{i0}=x^2F_{i2}$ and
$F_{i0}+F_{i2}z^2=(1+x^2z^2)F_{i0}$. Since
$B_{i0}+B_{i2}x=F_{i0}+F_{i2}z^2$, the sets $B_{i0}$ and $B_{i2}$
are unions of $2^{m-2}$ cosets of $\<x^2z^2\>$ and $A_{i0}$ and
$A_{i2}$ are unions of $2^{m-2}$ cosets of $\<x^2\>$ as $A_{is}$ is
the image of $B_{is}$ under the natural projection
$\overline{G}\rightarrow\overline{G}/\<z^2\>$. Since
$A_{i1}=\<y,F,x^2\>-A_{i0}$ and $A_{i3}=\<y,F,x^2\>-A_{i2}$,
$A_{i1}$ and $A_{i3}$ are also unions of $2^{m-2}$ cosets of
$\<x^2\>$. Clearly these cosets in $A_{i1}$ and $A_{i3}$ must come
from cosets of $\<x^2\>$ in $B_{i1}$ and $B_{i3}$ or cosets of
$\<x^2z^2\>$ in $B_{i1}$ and $B_{i3}$. Therefore $B_{i1}$ and
$B_{i3}$ are in unions of cosets of $\<x^2\>$ and cosets of
$\<x^2z^2\>$. If $B_{i1}$ and $B_{i3}$ contains no cosets of
$\<x^2z^2\>$, i.e. $B_{i1}+B_{i3}x=F_{i1}+F_{i3}z^2$ is the union of
$2^{m-1}$ cosets of $\<x^2\>$, let $\chi_1\in\widehat{\<x^2\>}$ such
that $\chi_1(x^2)=-1$ and $\chi_2\in\widehat{\<z^2\>}$ such that
$\chi_2(z^2)=-1$ and $\chi_1\chi_2$ and $\chi_2$ induces
$\C$-algebra homomorphisms
$$\chi_1\chi_2:\C[\<y,F,x^2,z^2\>]\rightarrow\C[\<y,F\>]$$
and $$\chi_2:\C[\<K,z^2\>]\rightarrow\C[K].$$ For every character
$\phi\in\widehat{\<y,F\>}$, one has
$\phi\chi_1\chi_2\in\widehat{\<y,F,x^2,z^2\>}$ which can be extended
to $\overline{G}$. Applying such an extension to $\overline{E}_i$,
one finds that $\phi\chi_1\chi_2(B_{i0})=\phi\chi_1\chi_2(B_{i2})=0$
for all $\phi$. Hence $\chi_1\chi_2(B_{i0})=\chi_1\chi_2(B_{i2})=0$
in $\C[\<y,F\>]$.  On the other hand, for
$\chi_0\in\widehat{\<x^2\>}$ with $\chi_0(x^2)=1$, one also has
$\chi_0\chi_2(B_{i0})=\chi_0\chi_2(B_{i2})=0$ as $B_{i0}$ and
$B_{i2}$ are unions of $2^{m-2}$ cosets of $\<x^2z^2\>$. Therefore
$\chi_2(B_{i0})=\chi_2(B_{i2})=0$ in $\C[\<y,F,x^2\>]$. This implies
that $\chi_2(F_{i0}+F_{i2}z^2)=0$ in $\C[K]$ or $F_{i0}=F_{i2}$ in
$K$. Since $F_{i0}$ and $F_{i2}$ are disjoint,
$F_{i0}=F_{i2}=\emptyset$, which contradicts the assumption that
$B_{i0}$ contains a coset of $\<x^2z^2\>$. Hence $B_{i1}$ or
$B_{i3}$ contains cosets of $\<x^2z^2\>$, i.e. $F_{i1}+F_{i3}z^2$
contains cosets of $\<x^2z^2\>$. Using arguments similar to those
for $F_{i0}+F_{i2}z^2$, one can show that $F_{i1}+F_{i3}z^2$ is the
union of $2^{m-1}$ cosets of $\<x^2z^2\>$. Therefore
$\overline{E}_i$ is the union of $2^m$ cosets of $\<x^2z^2\>$, i.e.
the coefficient of $x^2z^2$ in $\overline{E}_i\overline{E}^{(-1)}_i$
is $|\overline{E}_i|=2^{m+1}$. Thus for each $i\geqslant 2$, the
coefficient of $x^2z^2$ in $\overline{E}_i\overline{E}^{(-1)}_i$ is
either $0$ or $2^{m+1}$. Since the coefficient of $x^2z^2$ in
$\overline{E}_1\overline{E}^{(-1)}_1$ is $2^m$, the coefficient of
$x^2z^2$ in
$$\sum_{k=1}^{2^{m-1}}E_kE_k^{(-1)}$$ is congruent to $2^m$
$(\!\!\!\!\mod 2^{m+1})$. By Lemma~\ref{number}, one gets
$2^{2m-2}\equiv 2^m$ $(\!\!\!\!\mod 2^{m+1})$, a contradiction as
$m\geqslant 3$.

Therefor the group $N$ must be elementary abelian. By
Lemma~\ref{Eg}, Corollary~\ref{EI}, Lemma~\ref{noII} and
Lemma~\ref{I}, the group $G$ contains an index $2$ elementary
abelian subgroup $H\cong N\times N$. \qed

As we can see from the proof of Theorem~\ref{main}, the construction
of abelian $(2^{2m+1}(2^{m-1}+1), 2^m(2^m+1), 2^m)$-difference sets
for $m\geqslant 3$ is essentially unique and is given in
\cite{Dillon0} and \cite{Mc}.

\end{document}